\documentclass[11pt, reqno]{amsart}
\usepackage{multicol}
\usepackage[dvips]{graphicx}
\usepackage[ansinew]{inputenc}
\usepackage{amssymb}
\usepackage{amscd}
\usepackage{verbatim,ifthen}
\usepackage{color}
\usepackage{latexsym}
\usepackage{tikz}
\usepackage{tikz-cd}
\usepackage{mathrsfs}
\usepackage{wrapfig}
\usetikzlibrary{shapes}
\usepackage{color}
\usetikzlibrary{arrows.meta}
\usepackage{bbm}
\usetikzlibrary{matrix}
\usetikzlibrary{calc}
\usetikzlibrary{arrows,intersections}
\usepackage{pgfplots}
\usepackage{multicol}
\usepackage{array}
\newcolumntype{M}[1]{>{\centering\arraybackslash}m{#1}}

\usepackage[colorlinks, linkcolor=black, citecolor=magenta, linktocpage]{hyperref}
\addtocontents{toc}{\setcounter{tocdepth}{1}}

\usepackage{amsmath}

\numberwithin{equation}{section}

\let\oldtocsection=\tocsection
 
\let\oldtocsubsection=\tocsubsection

\renewcommand{\tocsection}[2]{\hspace{0em}\oldtocsection{#1}{#2}}
\renewcommand{\tocsubsection}[2]{\hspace{1em}\oldtocsubsection{#1}{#2}}

\usepackage{xcolor}

\def\XXint#1#2#3{{\setbox0=\hbox{$#1{#2#3}{\int}$ }
\vcenter{\hbox{$#2#3$ }}\kern-.6\wd0}}

\usepackage{eucal}
\usepackage{calc}  
\usepackage{enumitem} 
\usepackage{tensor}
\usepackage{graphicx,wrapfig,lipsum}
\usepackage{etoolbox}
\usepackage{marginnote}
\usepackage{lipsum}
\makeatletter
\patchcmd{\@mn@margintest}{\@tempswafalse}{\@tempswatrue}{}{}
\patchcmd{\@mn@margintest}{\@tempswafalse}{\@tempswatrue}{}{}
\reversemarginpar 
\makeatother

\usepackage{scrextend}
\newcommand*{\fullref}[1]{\hyperref[{#1}]{ \nameref*{#1} \autoref*{}}}
\makeatletter
\DeclareRobustCommand\widecheck[1]{{\mathpalette\@widecheck{#1}}}
\def\@widecheck#1#2{%
    \setbox\z@\hbox{\m@th$#1#2$}%
    \setbox\tw@\hbox{\m@th$#1%
       \widehat{%
          \vrule\@width\z@\@height\ht\z@
          \vrule\@height\z@\@width\wd\z@}$}%
    \dp\tw@-\ht\z@
    \@tempdima\ht\z@ \advance\@tempdima2\ht\tw@ \divide\@tempdima\thr@@
    \setbox\tw@\hbox{%
       \raise\@tempdima\hbox{\scalebox{1}[-1]{\lower\@tempdima\box
\tw@}}}%
    {\ooalign{\box\tw@ \cr \box\z@}}}
\makeatother
\title{The Schwarz lemma: An Odyssey}
\author{Kyle Broder}

\begin{document}

\maketitle

\begin{abstract}
Expository notes on the Schwarz lemma born out of some lectures given on the subject.
\end{abstract}

%\begin{multicols}{2}

\section{The birth of an empire}

The Schwarz lemma is one of the pioneering results in the function theory of one complex variable:

\subsection{Theorem}
A holomorphic map $f : \mathbb{D}(R_1) \to \mathbb{D}(R_2)$ fixing the origin, satisfies, for all $z \in \mathbb{D}(R_1)$, \begin{eqnarray}\label{Schwarz}
|f(z) | \ \leq \ \frac{R_2}{R_1} | z |.
\end{eqnarray}

The standard proof that we teach is the following: If $f(0)=0$, the function $g(z) := f(z)/z$ admits a holomorphic extension to all of $\mathbb{D}(R_1)$. Applying the maximum principle to $g(z)$ on each disk $| z | \leq R_1 - \varepsilon$, and letting $\varepsilon \to 0$ proves the statement. \\

Let us remark that this was not the proof originally given by Schwarz in \cite{SchwarzOriginal} (who proved the Schwarz lemma for one-to-one holomorphic maps). The proof here was first presented by Caratheodory \cite[p. 114, Note 13]{Caratheodory}, where it is attributed to Erhard Schmidt. 

In \eqref{Schwarz}, if we keep $R_2$ fixed, and let $R_1$ get arbitrarily large, we recover the following well-known corollary: 

\subsection{Corollary}
(Liouville's theorem). A bounded holomorphic function $f : \mathbb{C} \to \mathbb{C}$ assumes at most one value. \\

The Schwarz lemma is a local, finite statement about holomorphic maps. The Liouville theorem, in contrast, is a global statement obtained from letting $R_1 \to \infty$. This is the prototypical example of the so-called \textit{Bloch principle}; namely, the principle that any global statement concerning holomorphic maps arises from a stronger, finite version:
\begin{center}
\textit{Nihil est in infinito quod non prius fuerit in finito.}\footnote{There is nothing in the infinite that has not previously been in the finite. See \cite[p. 84]{BlochQuote}}
\end{center}

To further illustrate Bloch's principle, let us recall the following generalization of the Liouville theorem: The \textit{Picard theorem} states 

\subsection{Corollary}
A non-constant entire function $f : \mathbb{C} \to \mathbb{C}$ assumes all but possibly one value. \\

The corresponding finite version prophesized by the Bloch principle is the \textit{Schottky theorem}: 

\subsection{Theorem}
Let $f : \mathbb{D} \to \mathbb{C}$ be a holomorphic map which omits the values $0$ and $1$. Then $| f(z) |$ affords a bound in terms of $| f(0) |$ and $| z |$.\footnote{The original proof \cite{Schottky} gave no explicit bound for $|f(z)|$. }\\

The first instance of the Bloch principle was the following \textit{Valiron theorem}: 

\subsection{Corollary}
A non-constant entire function has holomorphic branches of the inverse in arbitrarily large Euclidean disks. \\

In \cite{Bloch}, Bloch improved Valiron's arguments and proved the underlying finite agent, which we now call the \textit{Bloch theorem}: 

\subsection{Theorem}
Every holomorphic function $f:\mathbb{D} \to \mathbb{C}$ has an inverse branch in some Euclidean disk of radius $\mathcal{B} | f'(0) |$, where $\mathcal{B} > 0$ is an absolute constant. 

\hfill

The constant $\mathcal{B}$ is now called \textit{Bloch's constant}, and its precise value remains unknown.

\section{The first divide: The Schwarz Pick Lemma}
As it stands, the Schwarz lemma controls the distortion of holomorphic maps between disks.  It was observed by Pick \cite{Pick1916} that the Schwarz lemma could be given a radically different interpretation.  Pick's first observation was that we do not require $f(0)=0$. Indeed, suppose $f : \mathbb{D} \to \mathbb{D}$ is a holomorphic self-map of the unit disk. For $\alpha \in \mathbb{D}$, the M\"obius transformation $$\varphi_{\alpha} : \mathbb{D} \to \mathbb{D}, \hspace*{1cm} \varphi_{\alpha}(z) := \frac{z-\alpha}{1-\overline{\alpha} z}$$ defines an automorphism of $\mathbb{D}$ which sends $\alpha$ the origin. The inverse of $\varphi_{\alpha}$ is, moreover, $\varphi_{\alpha}^{-1} = \varphi_{-\alpha}$. If $f(0) \neq 0$, we can produce a holomorphic self-map of $\mathbb{D}$ which fixes the origin by considering the composite map $$\varphi_{f(z)} \circ f \circ \varphi_{-z}.$$ Setting $w = \varphi_{-z}(\zeta)$, the familar Schwarz lemma gives \begin{eqnarray*}
| \varphi_{f(z)} \circ f(w) | & \leq & | \varphi_z(w) |.
\end{eqnarray*}

Explicitly, this reads \begin{eqnarray*}
\left| \frac{f(w) - f(z)}{1 - \overline{f(z)} f(w) } \right| & \leq & \left| \frac{w-z}{1-\overline{z}w} \right|.
\end{eqnarray*}

The function $$d_{\text{H}} : \mathbb{D} \times \mathbb{D} \to \mathbb{R}, \hspace*{1cm} d_{\text{H}} (z,w) : = \left| \frac{z-w}{1 - \overline{w} z}  \right|$$ defines a distance function on $\mathbb{D}$, the \textit{pseudo-hyperbolic distance}. The function $d_{\text{H}}$ defines an honest distance function\footnote{i.e., a symmetric non-degenerate function satisfying the triangle inequality. }.  It does not, however, come from integrating a Riemannian metric.

\subsection{Reminder: Distance functions from Riemannian metrics}
A Riemannian metric $g$ on a smooth manifold $M$ is a smoothly varying family of positive-definite quadratic forms $g_p$ on each tangent space $T_p M$. The metric permits us to calculate the length of smooth curves $\gamma : [0,1] \to M$ in $M$ by integrating their tangent vectors: \begin{eqnarray}\label{Length}
\text{length}_g(\gamma) : =   \int_0^1 | \gamma'(t) |_{g,\gamma(t)} dt.
\end{eqnarray}

We can subsequently construct a distance function $\text{dist}_g : M \times M \to \mathbb{R}$ from the formula $$\text{dist}_g(p,q) \ := \ \inf_{\gamma} \int_0^1 | \gamma'(t) |_{g,\gamma(t)} dt.$$ Here, the infimum is taken over all smooth curves $\gamma : [0,1] \to M$ such that $\gamma(0) =p$ and $\gamma(1)=q$.

\subsection{The Poincar\'e metric}
On the unit disk $\mathbb{D} = \{ | z | < 1 \}$, we have the \textit{Poincar\'e metric}: $$\rho \ : = \ \frac{| dz |}{(1- | z |^2)}.$$

Let us compute the length of a line segment $\ell$ connecting $0$ with $w \in (0,1)$, with respect to the Poincar\'e metric. Indeed, parametrizing $\ell$ by $\gamma(t) = t w$, we see that \begin{eqnarray*}
\text{length}_{\rho}(\ell) \ = \  \int_{\ell} d \rho &=& \int_0^1 \frac{w}{1- w^2 t^2} dt \\ 
&=&  \frac{1}{2} \log \frac{1+w}{1-w} \ = \ \tan^{-1}(w).
\end{eqnarray*}

It is an elementary exercise to show that the associated \textit{Poincar\'e distance function} is given by $$\text{dist}_{\rho}(z,w) \ = \ \tan^{-1} \left( \frac{z-w}{1-z \overline{w}} \right) .$$

Summarizing this discussion,  and replacing the pseudo-hyperbolic distance with the Poincar\'e distance, we have recovered the theorem of Pick \cite{Pick1916}:

\subsection{Theorem}
(Schwarz--Pick lemma). Let $f : \mathbb{D} \to \mathbb{D}$ be a holomorphic self-map of the unit disk. Then for all $z,w \in \mathbb{D}$, \begin{eqnarray*}
\text{dist}_{\rho}(f(z), f(w)) & \leq & \text{dist}_{\rho}(z,w).
\end{eqnarray*}

That is, with respect to the Poincar\'e distance, all holomorphic maps are distance-decreasing.

The Schwarz--Pick lemma holds equally well with the Poincar\'e distance function replaced by the pseudo-hyperbolic distance $d_{\text{H}}$. It will be clear in subsequent developments why we have chosen to state the Schwarz--Pick lemma in terms of the Poincar\'e distance. \footnote{Independent of this, however, the statement given by Pick in \cite{Pick1916} was in terms of the Poincar\'e metric, not the pseudo-hyperbolic distance.}

\section{Aside: Intrinsic metrics}

The Schwarz--Pick lemma casts the focus of the Schwarz lemma away from holomorphic maps. Instead, the Poincar\'e distance is the object of focus. There are many examples of such (pseudo-)distance functions, generalizing the Poincar\'e metric:

The \textit{Carath\'eodory metric} $d_{\Omega}$ is the pseudo-distance function defined on a domain $\Omega \subseteq \mathbb{C}$ by $$d_{\Omega}(z,w) := \sup_f \text{dist}_{\rho}(f(z),f(w)),$$ where, the supremum is taken over all holomorphic maps $f : \Omega \to \mathbb{D}$. Clearly, all holomorphic maps are distance-decreasing wth respect $d_{\Omega}$. The Carath\'eodory pseudo-distance does not define a distance function, in general, since it may be degenerate (e.g., when $\Omega = \mathbb{C}$). 

The most notable example of a pseudo-distance for which all holomorphic maps are distance-decreasing, and which can be defined on any complex manifold $M$,  is the \textit{Kobayashi pseudo-distance}: \begin{eqnarray*}
d_K(p,q) & : = & \inf \sum_{j=1}^m \text{dist}_{\rho}(s_j,t_j),
\end{eqnarray*}

where the infimum is taken over all $m \in \mathbb{N}$, all pairs of points $\{ s_j, t_j \}$ in $\mathbb{D}$, and all collections of holomorphic maps $f_j : \mathbb{D} \to M$ such that $f_1(s_1) = p$, $f_m(t_m) = q$, and $f_j(t_j) = f_{j+1}(s_{j+1})$.

A complex manifold is said to be \textit{Kobayashi hyperbolic} if the Kobayashi pseudo-distance is an honest distance function, i.e., non-degenerate.

A compact Riemann surface $\Sigma_g$ of genus $g \geq 2$ is Kobayashi hyperbolic. Note that in this case, the Carath\'eodory pseudo-distance vanishes identically!

\section{The second rift: The Ahlfors Schwarz lemma}

In contrast with the pseudo-hyperbolic distance, the Poincar\'e distance comes from a Riemannian metric. In particular, one can use the machinery of differential geometry to study holomorphic maps. From this viewpoint, the Poincar\'e metric has some very special properties:

\subsection{Reminder: Gauss Curvature}
Let $g = \lambda | dz |$ be a Hermitian metric on a domain $\Omega \subseteq \mathbb{C}$, or more generally, a Riemann surface. The \textit{Gauss curvature} of $g$ is the function $$K_g \ = \ - \frac{1}{\lambda^2} \Delta \log \lambda \ = \ - \frac{1}{\lambda^2} \frac{\partial^2}{\partial z \partial \overline{z}} \log \lambda.$$ 

In particular, if $K_g \leq -C$, then setting $u=\log(\lambda)$, we have $\Delta u \geq C e^{2u}$.\\

Applying the above formula, the curvature of Poincar\'e metric $$\rho = \frac{| dz |}{(1- | z |^2)}$$ is seen to be $K_{\rho} \equiv -4$. That is, the Poincar\'e metric has constant (negative) Gauss curvature.

If we suppose that $g = e^u | dz |$ is some metric of negative curvature, in the sense that $K_g \leq -4$, it is natural to ask how $g$ compares with the Poincar\'e metric. If we let $R = 1 - \varepsilon$ for some $\varepsilon \in (0,1)$, the Poincar\'e metric on $\mathbb{D}(R)$ is given by $\rho_R = e^v | d z | :=   \frac{R}{R^2-| z |^2} | dz |$. Since $K_{\rho_R} \equiv -4$, we have $\Delta v = 4 e^{2v}$. Hence, on the open set $\Omega := \{ z \in \mathbb{D}(R) : u(z) > v(z) \}$, the function $u-v$ is subharmonic: $\Delta(u - v) \geq e^{2u} - e^{2v}$. By the maximum principle, $u-v$ cannot achieve an interior maximum, and hence, the supremum must be approached on the boundary. But $\Omega$ cannot have boundary points on $| z | = R$, since $v \to \infty$ as $| z | \to R$. By continuity, at a boundary point $z$ of $\Omega$ with $| z | < R$,  we have $u-v=0$, yielding a contradiction if $\Omega$ is non-empty. We therefore deduce that $u(z) \leq v(z)$ for all $| z | < R$. Letting $\varepsilon \to 0$ recovers the following theorem of Ahlfors:

\subsection{Theorem}
(Ahlfors--Schwarz lemma). Let $(\Sigma, g)$ be a Riemann surface with Gauss curvature $K_g \leq -4$.  Then for all holomorphic maps $f : \mathbb{D} \to \Sigma$,  $$\text{dist}_g(f(z),f(w)) \ \leq \ \text{dist}_{\rho}(z,w),$$ for all $z,w \in \mathbb{D}$.\\

In the collected works of Ahlfors \cite[p. 341]{AhlforsCollectedWorks}, one finds the following reflection concerning his version of the Schwarz lemma: Ahlfors confesses that his generalization of the Schwarz lemma had ``\textit{more substance that I was aware of}", but ``\textit{without applications, my lemma would have been too lightweight for publication}". The applications that Ahlfors alludes to here were the following: First, a proof of Schottky's theorem with definite numerical bounds: If $f : \mathbb{D} \to \mathbb{C}$ is a holomorphic function such that $f(\mathbb{D}) \cap \{ 0,1 \} = \emptyset$, then \begin{eqnarray*}
\log | f(z) | & < & \frac{1+\vartheta}{1 - \vartheta} \left( 7 + \max \{ 0, \log | f(0) | \} \right),
\end{eqnarray*}

for all $ | z | \leq \vartheta < 1$. The second was an improved lower bound on the Bloch constant: $$\mathcal{B} \ \geq \ \frac{\sqrt{3}}{4}.$$

\section{The Aftermath of the Ahlfors Schwarz Lemma}

The Ahlfors incarnation again redirects the focus of the Schwarz lemma away from holomorphic maps, and away from the distance function; it is the curvature that is the primary object of focus.

\subsection{Reminder: Curvature}
Let $(M^n, g)$ be a Riemannian manifold of dimension $n$. A \textit{connection} (or \textit{covariant derivative}) $\nabla$ on $M$ provides a coordinate-free tool for computing the directional derivative of a vector field on $M$.

A connection is said to be \textit{metric} if $\nabla g =0$, i.e., $w(g(u,v)) = g(\nabla_w u, v) + g(u, \nabla_w v)$ for all vector fields $u,v,w$. The \textit{Levi-Civita connection} is the unique metric connection on $TM$ which is torsion-free, i.e., $T(u,v) : = \nabla_u v - \nabla_v u - [u,v]=0$, where $[ \cdot, \cdot ]$ denotes the Lie bracket.

The obstruction to covariant derivatives commuting is measured by the \textit{Riemannian curvature}: $$\text{Rm}_g(u,v) w := \nabla_u \nabla_v w - \nabla_v \nabla_u w - \nabla_{[u,v]} w,$$ where $\nabla$ is the Levi-Civita connection. Using the Riemannian metric $g$, we can write $\text{Rm}$ as $$\text{Rm}_g(u,v,w,z) : = g(\text{Rm}_g(u,v)w,z),$$ rendering it scalar-valued. If $e_1, ..., e_n$ denote a local frame for $TM$ near a point $p$, then we can write $$\text{Rm}_g(u,v,w,z) = \sum R_{i j k \ell} u_i v_j w_k z_{\ell},$$ 

where $R_{ijk\ell}$ are the components of $\text{Rm}_g$ with respect to the frame.

If $u,v \in T_p M$ are vectors which form an orthonormal basis for a two-dimensional subspace of $T_p M$, we define the \textit{sectional curvature} by $$\text{Sec}_g(u,v) := \text{Rm}_g(u,v,v,u).$$ 

The metric trace of the curvature tensor yields the \textit{Ricci curvature} $$\text{Ric}_g(u,v) = \sum_i \text{Rm}_g(u,e_i, v, e_i).$$ With respect to a frame, the components of $\text{Ric}_g$ are $$\text{Ric}_{k\ell} = \sum_{i,j} g^{ij}R_{ki\ell j}.$$

\subsection{Harmonic maps}
Let $f : (M,g) \to (N, h)$ be a smooth map of Riemannian manifolds. The \textit{energy} of $f$ is defined $$\mathcal{E}(f) := \int_M | d f |^2 dV_g.$$ The Euler--Lagrange operator associated to $\mathcal{E}$ is the \textit{tension field} of $f$, namely, $\tau(f) = \text{div}(d f)$. We say that $f$ is \textit{harmonic} if $\tau(f)=0$. In \cite{EellsSampson}, it was shown that if $f$ is harmonic, then \begin{eqnarray}\label{HarmonicMapBochner}
\Delta | df |^2 &=& | \nabla d f |^2 + \langle \text{Ric}_g(\nabla_v f), \nabla_v f \rangle \\
&& \hspace*{3cm} - \langle \text{Riem}_h(\nabla_v f, \nabla_w f) \nabla_v f, \nabla_w f \rangle.\nonumber 
\end{eqnarray}

In particular, if $(M, g)$ is compact with $\text{Ric}_g > 0$, and $\text{Sec}_h \leq 0$, then $f$ is constant.

\subsection{Curvature in the presence of a complex structure}
For harmonic maps $f : (M, g) \to (N,h)$ between Riemannian manifolds, we required a lower bound on the Ricci curvature of the source metric $g$ and an upper bound on the sectional curvature of the target metric $h$. It is hard to do much better than this, in general. In the complex analytic setting, however, one can say a lot more:

\subsection{Complex manifolds}
A complex manifold is a smooth manifold $M$ with local holomorphic coordinates in a neighborhood of each point. By a theorem of Newlander--Nirenberg, this is equivalent to existence of an endomorphism $J : TM \to TM$ satisfying $J^2 = - \text{id}$ together with an additional integrability criterion: The relation $J^2 = - \text{id}$ splits the (complexified) tangent bundle into a direct sum of eigenbundles $T^{\mathbb{C}} M = T^{1,0}M \oplus T^{0,1}M$, corresponding to the eigenvalues $\pm \sqrt{-1}$. The complex structure $J$ is said to be \textit{integrable} if $T^{0,1}M$ is closed under the Lie bracket.

In contrast with the Riemannian theory, the Levi-Civita connection is not the natural connection on the tangent bundle of a complex manifold. Indeed, we say that a metric connection $\nabla$ is Hermitian if $\nabla J=0$. The preferred connection here is the \textit{Chern connection}: the unique Hermitian connection whose torsion $T$ satisfies $T(Ju,v) = T(u,Jv)$ for all $u,v \in TM$. A complex manifold, therefore, has two distinguished connections on its tangent bundle. The class of manifolds for which these two connections coincide are the well-known \textit{K\"ahler manifolds}. Such manifolds exist in an almost-baffling abundance: Complex projective space $\mathbb{P}^n$, Euclidean space $\mathbb{C}^n$, the complex ball $\mathbb{B}^n$ all admit K\"ahler structures. Complex submanifolds of K\"ahler manifolds remain K\"ahler. Hence, all projective manifolds and Stein manifolds are K\"ahler.

\subsection{Curvature of Hermitian Metrics}
For a Hermitian manifold $(M, g)$, the \textit{holomorphic bisectional curvature} is defined $$\text{HBC}_g(u,v) :=  R(u, \overline{u}, v, \overline{v}) =   \sum R_{i \overline{j} k \overline{\ell}} u_i \overline{u}_j v_k \overline{v}_{\ell},$$ for unit $(1,0)$--tangent vectors $u = (u_1, ..., u_n)$, $v=(v_1, ..., v_n) \in T_p^{1,0}M$, $p \in M$,  where $R$ denotes the curvature of the Chern connection.  The terminology introduced by Goldberg--Kobayashi \cite{GoldbergKobayashiHBC} is justified by the fact that if $g$ is K\"ahler, then the holomorphic bisectional curvature can be written as a sum of two sectional curvatures. 

The bisectional curvature is very restrictive on the geometry. For instance, by the Mori \cite{Mori}, Siu--Yau \cite{SiuYau} solution of the Frankel conjecture: if $(X,\omega)$ is a compact K\"ahler manifold with a K\"ahler metric of positive bisectional curvature then $X$ is biholomorphic to $\mathbb{P}^n$.  If, in addition, the K\"ahler metric is Einstein,  by Goldberg--Kobayashi \cite{GoldbergKobayashiHBC}, the metric is biholomorphically isometric to $\mathbb{P}^n$ with the Fubini--Study metric.  

The restriction of the holomorphic bisectional curvature to the diagonal yields the \textit{holomorphic sectional curvature}: $$\text{HSC}_g(v) : = R(v, \overline{v}, v, \overline{v}) = \sum R_{i \overline{j} k \overline{\ell}} v_i \overline{v}_j v_k \overline{v}_{\ell},$$ for a unit $(1,0)$--tangent vector $v \in T_p^{1,0}M$, $p \in M$.

The holomorphic sectional curvature is comparatively much weaker, but does have important implications on the complex geometry: By the Ahlfors--Schwarz lemma,  a complex manifold with a Hermitian metric of negative holomorphic sectional curvature (bounded above by a negative constant) is Brody hyperbolic (every entire curve, i.e., holomorphic map $\mathbb{C} \to X$, is constant).  If, in addition, $X$ is compact,  then $X$ is Kobayashi hyperbolic. On the other hand,  by a theorem of Yang \cite{YangRC}, a compact K\"ahler manifold with a K\"ahler metric of positive holomorphic sectional curvature is rationally connected (any two points lie in the image of a rational curve, i.e., holomorphic map $\mathbb{P}^1 \to X$).

\subsection{Higher-dimensional Schwarz Lemmas}
The first general result in higher-dimensions was obtained by Chern \cite{Chern66} and Lu \cite{Lu}. Let $f : (M,g) \to (N,h)$ be a holomorphic map between Hermitian manifolds.  Choose coordinates $(z_1, ..., z_m)$ at a point $p \in M$ and $(w_1, ..., w_n)$ at $f(p) \in N$.  If $f = (f^1, ..., f^n)$, we write $f_i^{\alpha} := \frac{\partial f^{\alpha}}{\partial z_i}$.  With this notation maintained,  we recall the following general formula for the Laplacian of $| \partial f |^2$ obtainted by Chern \cite{Chern66} and Lu \cite{Lu}:

\subsection{Theorem}
(\cite{Chern66}, \cite[Theorem 4.1]{Lu}). Let $f: (M, g) \to (N, h)$ be a holomorphic map between Hermitian manifolds.  Then \begin{eqnarray}\label{LuFormula}
\frac{1}{2} \Delta_g | \partial f |^2 &=& | \nabla \partial f |^2 + g^{i \overline{j}} R_{i \overline{j} k \overline{\ell}}^g g^{k \overline{q}} g^{p \overline{\ell}} h_{\alpha \overline{\beta}} f_p^{\alpha} \overline{f_q^{\beta}} \\
&& \hspace*{3cm} - R_{\alpha \overline{\beta} \gamma \overline{\delta}}^h \left( g^{i \overline{j}} f_i^{\alpha} \overline{f_j^{\beta}} \right) \left( g^{p \overline{q}} f_p^{\gamma} \overline{f_q^{\gamma}} \right).\nonumber 
\end{eqnarray}

The right-hand side of formula \eqref{LuFormula} has three terms: $| \nabla \partial f |^2$ is the pointwise norm squared of the \textit{second fundamental form of $f$}.\footnote{This terminology is justified by the fact that if $f$ is a harmonic immersion between Riemannian manifolds, then $| \nabla \partial f |^2$ is the pointwise norm squared of the familiar second fundamental form.} The second term is controlled by assuming a lower bound on the second Chern--Ricci curvature\footnote{In the non-K\"ahler case, there are three distinct Ricci curvatures: $\text{Ric}^{(1)}_{i \overline{j}}  : =  g^{k \overline{\ell}} R_{i \overline{j} k \overline{\ell}}$, $\text{Ric}^{(2)}_{k \overline{\ell}}  : =  g^{i \overline{j}} R_{i \overline{j} k \overline{\ell}}$, and $\text{Ric}^{(3)}_{k \overline{j}}  : =  g^{i \overline{\ell}} R_{i \overline{j} k \overline{\ell}}.$} of $g$. Indeed, if we suppose that $\text{Ric}_g^{(2)} \geq - C g$ for some constant $C>0$, then since $g^{i \overline{j}} R_{i \overline{j}k\overline{\ell}}^g = \text{Ric}_{k \overline{\ell}}^{(2)}$, we see that \begin{eqnarray}\label{SecondChernRicciLowerBound}
g^{i \overline{j}} R_{i \overline{j}k\overline{\ell}} g^{k \overline{q}} g^{p \overline{\ell}} h_{\alpha \overline{\beta}} f_p^{\alpha} \overline{f_q^{\beta}} \geq - C g_{k \overline{\ell}} g^{k \overline{q}} g^{p \overline{\ell}} h_{\alpha \overline{\beta}} f_p^{\alpha} \overline{f_q^{\beta}} = -C g^{p \overline{q}} h_{\alpha \overline{\beta}} f_p^{\alpha} \overline{f_q^{\beta}} = -C | \partial f |^2.
\end{eqnarray}

The last term in formula \eqref{LuFormula} has caused the most confusion and is indeed the most mysterious. Lu \cite{Lu} falsely identifies this term as the holomorphic sectional curvature. In the K\"ahler setting, Royden \cite{Royden} showed that this term is controlled if one has a non-positive upper bound on the holomorphic sectional curvature of $h$.  The crux of Royden's argument is the following polarization result, which is commonly referred to as \textit{Royden's trick} in the literature (see, e.g., \cite{YangZhengRBC,WuYau1,TosattiYang}):

\subsection{Proposition}
Let $\xi_1, ..., \xi_{\nu}$ be othogonal tangent vectors.  Let $S(\xi, \overline{\eta}, \zeta, \overline{\omega})$ be a symmetric bi-Hermitian form in the sense that \begin{eqnarray*}
S(\xi, \overline{\eta}, \zeta, \overline{\omega}) \ = \ S(\zeta, \overline{\eta}, \xi, \overline{\omega}), \hspace*{1cm} \text{and} \hspace*{1cm} S(\eta, \overline{\xi}, \omega, \overline{\zeta}) \ = \ \overline{S}(\xi, \overline{\eta}, \zeta, \overline{\omega}).
\end{eqnarray*}

If $S(\xi, \overline{\xi}, \xi, \overline{\xi}) \leq \kappa \| \xi \|^4$, then $$\sum_{\alpha, \beta} S(\xi_{\alpha}, \overline{\xi}_{\alpha}, \xi_{\beta}, \overline{\xi}_{\beta}) \ \leq \ \frac{1}{2} \kappa \left[ \left( \sum_{\alpha} \| \xi_{\alpha} \|^2 \right)^2 + \sum_{\alpha} \| \xi_{\alpha} \|^4 \right].$$ Further, if $\kappa \leq 0$, then \begin{eqnarray}\label{Roydennonpositive}
\sum_{\alpha, \beta} S(\xi_{\alpha},\overline{\xi}_{\alpha}, \xi_{\beta}, \overline{\xi}_{\beta}) & \leq & \kappa \frac{n+1}{2n} \left( \sum_{\alpha} \| \xi_{\alpha} \|^2 \right)^2.
\end{eqnarray}

Applying this polarization argument to the target curvature term, assuming the metric is K\"ahler, shows that the target curvature term can indeed be controlled by the holomorphic sectional curvature.  It is worth noting that the same argument holds for a larger class of metrics, namely metrics which are (Chern) \textit{K\"ahler-like} \cite{YangZhengCurvature}: A Hermitian metric is said to be (Chern) K\"ahler-like if its Chern curvature tensor has the symmetries of the K\"ahler curvature tensor. The only known example of a non-K\"ahler (Chern) K\"ahler-like metric on a compact Hermitian manifold is the Iwasawa threefold, whose Chern curvature tensor vanishes identically.  Hence, for trivial reasons,  the Iwasawa threefold is (Chern) K\"ahler-like. It is expected, however, that such metrics should exist in abundance.\footnote{I would like to thank Professor Fangyang Zheng for many valuable communications on this topic.} For a general Hermitian metric, however, Royden's technique cannot be used, and it is not clear whether the holomorphic sectional curvature gives suitable control.

Combining \eqref{SecondChernRicciLowerBound} and \eqref{Roydennonpositive}, we see that if $g$ is a K\"ahler metric with $\text{Ric}_g \geq -C$ and $h$ is a K\"ahler metric with $\text{HSC}_h \leq - \kappa$, then \begin{eqnarray}\label{RoydenBochCalc}
\frac{1}{2} \Delta_g | \partial f |^2 &=& | \nabla \partial f |^2 - C | \partial f |^2 + \kappa \frac{r+1}{r} | \partial f |^4 \ \geq \ - C | \partial f |^2 + \kappa \frac{r+1}{r} | \partial f |^4, 
\end{eqnarray}

where $r$ is the rank of $\partial f$.

To deduce an estimate on $| \partial f |^2$ from \eqref{RoydenBochCalc}, we require the maximum principle. If $M$ is compact, then we have access to the maximum principle. In particular, applying \eqref{RoydenBochCalc}, we discover:

\subsection{Theorem}
(\cite{Chern66, Lu, Royden}). Let $f : (M, g) \to (N, h)$ be a holomorphic map of rank $r$ between K\"ahler manifolds. Assume $\text{Ric}_{g} \geq - C$ and $\text{HSC}_{h} \leq - \kappa$ for some constants $C >0$ and $\kappa \geq 0$. If $M$ is compact, then \begin{eqnarray*}
| \partial f |^2 & \leq & \frac{2C r}{(r+1) \kappa}.
\end{eqnarray*}

In particular, there are non-constant holomorphic maps $f : (M,g) \to (N,h)$ from a compact K\"ahler manifold with $\text{Ric}_g>0$ (or compact Hermitian with $\text{Ric}_g^{(2)}>0$) into a K\"ahler manifold with $\text{HSC}_h <0$.

At this point, let us emphasize that although we motivated the Schwarz lemma almost exclusively from the point of view of complex function theory (i.e., understanding estimates on the derivative of holomorphic maps), the Schwarz lemma is an essential tool (and often the bottleneck) to a number of problems in complex differential geometry.  One particularly important instance is the Wu--Yau theorem \cite{Wong, HeierLuWong, WongWuYau, WuYau1, WuYau2, DiverioTrapani, TosattiYang, YangZhengRBC, BroderSBC, DiverioSurvey}: 

\subsection{Theorem}
Let $(M, g)$ be a compact K\"ahler manifold with $\text{HSC}_g \leq - \kappa <0$ for some constant $\kappa>0$.  Then the canonical bundle $K_M$ is ample. In particular,  $M$ is projective and canonically polarized.  \\

By the Aubin--Yau theorem \cite{Aubin, Yau1976}, it follows that a compact K\"ahler manifold with negative holomorphic sectional curvature has a K\"ahler--Einstein metric of negative Ricci curvature.  The crux of the argument is to produce a K\"ahler--Einstein metric as a solution to a certain Monge--Amp\`ere equation (see \cite{WuYau1} for details) via a limiting process.  If $\varepsilon$ is the parameter in which a limit is taken,  then the argument hinges upon obtaining uniform estimates on all derivatives of the metric (independent of $\varepsilon$).  If one has a uniform $\mathcal{C}^2$--estimate (which is given exactly by the Schwarz lemma), then by standard elliptic theory (see \cite{TosattiYang} for a nice treatment of this part), all the higher-order estimates are obtained.  Hence, the bottleneck to obtaining higher-order a priori estimates for complex Monge--Amp\`ere equations is typically the Schwarz lemma. We invite the reader to consult \cite{DiverioSurvey, WuYau1, TosattiYang} for more details on the structure of the argument. \\

Although the pointwise equality \eqref{LuFormula} does not require the source manifold to be compact, in the absense of compactness, there is no guarentee that a maximum exists, and the maximum principle cannot be applied directly. One way to circument this in this in the non-compact case, is to consider manifolds with a certain exhaustion property: A manifold $M^n$ has the \textit{$K$--exhaustion property} \cite[$\S 5$]{Lu} if $M$ is exhausted by a sequence of open submanifolds $M_1 \subset M_2 \subset \cdots \subset M$, with compact closures and such that: (i) for each $k \in \mathbb{N}$, there is a smooth function $v_k \geq 0$ on $M_k$ with $\frac{1}{2} \Delta v_k  \leq  \frac{R}{n} + K \exp(v_k)$, for some fixed constant $K>0$; (ii) if $p_i$ is a divergent sequence\footnote{An infinite sequence $p_i$ in $M_k$ is said to be \textit{divergent} if every compact open set in $M_k$ contains only a number of points in this sequence.} of points in $M_k$, then $v_k(p_i) \to \infty$. The unit ball $\mathbb{B}^n$ has the $K$--exhaustion property with $K = 2n(n+1)$.  As a consequence, we can apply \eqref{RoydenBochCalc} if the source manifold is $\mathbb{B}^n$:

\subsection{Theorem}\label{CLSchwarzEstimate}
(\cite{Chern66}, \cite{Lu}). Let $f : \mathbb{B}^n \to N$ be a holomorphic map, where $\mathbb{B}^n$ is equipped with the standard metric, and $N$ is a K\"ahler manifold. If $N$ is equipped with a K\"ahler metric of negative holomorphic sectional curvature, then $f$ is distance-decreasing. \\

It is worth emphasizing that this is exactly why Chern considers this rather restrictive case in \cite{Chern66}. The $K$--exhaustion property is, of course, very restrictive. The breakthrough that was required for the Schwarz lemma in the non-compact case was made by Omori \cite{Omori} and Yau \cite{YauMaximumPrinciple}:

\subsection{Theorem}\label{YauMaximum}
Let $(M,g)$ be a complete Riemannian manifold. Assume $\text{Ric}_g \geq - C$ for some $C \in \mathbb{R}$. Let $f : M \to \mathbb{R}$ be a smooth function which is bounded above. Then for any $\varepsilon>0$, there is a point $p \in M$ such that $| f(p) - \sup_M f | < \varepsilon$,  $\| \text{grad} \vert_p f \| < \varepsilon$,  and $\Delta \vert_p f < \varepsilon$. \\

\iffalse

The subject of maximum principles for non-compact Riemannian manifolds is rich: Omori \cite{Omori} was the first to show that the naive extensions of the maximum principle to complete Riemannian manifolds did not hold. Omori established the above theorem under the more restrictive assumption of a lower bound on the sectional curvature. The lower bound on the Ricci curvature given in Yau's formulation \cite{YauMaximumPrinciple} is, of course, a substantial improvement. There are results due to \cite{Borbely, ChenXin, KimLee,RattiRigoliSetti}, however, indicating that one needs only control the rate of descent of the Ricci curvature as it approaches $-\infty$.

\fi

With the maximum principle now available in this level of generality, we have:

\subsection{Theorem}
(\cite{YauSchwarzLemma}). Let $f : (M,g) \to (N,h)$ be a holomorphic map from a complete K\"ahler manifold wth $\text{Ric}_g \geq - A$, to a Hermitian manifold $(N, h)$ with $\text{HBC}_h \leq - B <0$. Then $$| \partial f |^2 \ \leq  \ \frac{A}{B}.$$ 

\hfill

It should be noted that Yau's result appears before Royden's Schwarz lemma \cite{Royden}. On the other hand,  Yau remarks on page 1 of \cite{YauSchwarzLemma} that Royden could already improve the bisectional curvature bound to a holomorphic sectional curvature bound if the target metric is K\"ahler.

\subsection{The Aubin--Yau Inequality}
All forms of the Schwarz lemma so far have arisen (more or less) from applying the maximum principle to \eqref{LuFormula} or \eqref{RoydenBochCalc}. If we additionally assume that $f$ is biholomorphic, however, we have another Laplacian at our disposal; namely, the target metric Laplacian: $\Delta_{h} = \text{tr}_{h} \sqrt{-1} \partial \overline{\partial}$. The Schwarz lemma with the target metric Laplacian was first considered in \cite{Aubin, Yau1976}, and hence is referred to as the Aubin--Yau second-order estimate:

\subsection{Theorem}
(\cite{Aubin,Yau1976}). Let $f : (M^n, g) \longrightarrow (N, h)$ be a holomorphic map, which is biholomorphic onto its image. Assume $\text{HBC}_{g} \geq - \kappa$ and $\text{Ric}_{h}^{(2)} \leq - C_1 h + C_2 (f^{-1})^{\ast} g$ for constants $\kappa, C_1,C_2$, with $C_1 >0$. Then, if $M$ is compact,  \begin{eqnarray*}
| \partial f |^2 & \leq & \frac{nC_2 + \kappa}{C_1}.
\end{eqnarray*}

\subsection{A Unified Approach}
In \cite{RubinsteinKE}, Rubinstein gave a unified treatment of the Chern--Lu and Aubin--Yau inequalities in the K\"ahler setting. The underlying philosophy is that these theorems are best understood via holomorphic maps, and not as abstract tensor calculations. This paradigm extends to the study of holomorphic maps between Hermitian manifolds and was made more transparent in \cite{BroderSBC} using the ideas in \cite{YangZhengRBC}. In effect, this goes back to the original formulation of the Schwarz lemma: the focus should be placed on the holomorphic map.\footnote{As an aside, we must mention that the use of the Chern--Lu inequality for $\mathcal{C}^2$--estimates came to prominence only after the work of Jeffres--Mazzeo--Rubinstein \cite{JMR2016} (c.f., \cite[Theorem 2.2]{CDS1}). Prior to \cite{JMR2016},  the Aubin--Yau inequality was the primary tool.}

The Bochner formula details how to compute the complex Hessian of (the pointwise norm squared of) a holomorphic section $\sigma \in H^0(\mathcal{E})$ of a holomorphic vector bundle $\mathcal{E}$: \begin{eqnarray}\label{BOCHBASIC}
\sqrt{-1} \partial \overline{\partial} | \sigma |^2 &=& | \nabla \sigma |^2 -  \sqrt{-1} \langle R^{\mathcal{E}} \sigma, \sigma \rangle,
\end{eqnarray}

where $R^{\mathcal{E}}$ is the curvature of the Chern connection $\nabla$ on $\mathcal{E}$. When $\sigma = \partial f$, and $\mathcal{E} = T^{\ast} M \otimes f^{\ast} T^{1,0}N$, we have \begin{eqnarray}\label{BochnerMotivation}
\sqrt{-1} \partial \overline{\partial} | \partial f |^2 & = & \langle \nabla^{1,0} \partial f, \nabla^{1,0} \partial f \rangle  - \sqrt{-1} \langle R^{T^{\ast}M \otimes f^{\ast} T^{1,0}N} \partial f, \partial f \rangle.
\end{eqnarray}

Since the curvature of the tensor product of bundles splits additively, we get opposing contributions to the curvature from the source and target metrics: \begin{eqnarray}\label{SplitMotivation}
R^{T^{\ast}M \otimes f^{\ast} T^{1,0}N} &=& - R^{T^{1,0}M} \otimes \text{id} + \text{id} \otimes f^{\ast} R^{T^{1,0}N}.
\end{eqnarray}

Taking the trace of \eqref{BochnerMotivation} with respect to $g$, we recover the Chern--Lu formula \eqref{LuFormula}. As we indicated previously, we want to understand the target curvature term: \begin{eqnarray}\label{MOTIVATIONCURVATURE}
g^{i \overline{j}} g^{p \overline{q}} R_{ \alpha \overline{\beta}\gamma \overline{\delta}}^h f_i^{\alpha} \overline{f_{j}^{\beta}} f_p^{\gamma} \overline{f_{q}^{\delta}}.
\end{eqnarray}

Choose coordinates $(z_1, ..., z_n)$ centered at a point $p \in M$ and $(w_1, ..., w_n)$ at $f(p) \in N$ such that $g_{i \overline{j}} =\delta_{ij}$ and $h_{\alpha \overline{\beta}} = \delta_{\alpha\beta}$ at $p$ and $f(p)$, respectively. If $f = (f^1, ..., f^n)$, then with $f_i^{\alpha} = \frac{\partial f^{\alpha} }{\partial z_i}$, the coordinates can be chosen such that $f_i^{\alpha} = \lambda_i \delta_i^{\alpha}$, where $\lambda_1 \geq \lambda_2 \geq \cdots \geq \lambda_r \geq \lambda_{r+1} = \cdots = 0$, and $r$ is the rank of $\partial f = (f_i^{\alpha})$. Hence, \eqref{MOTIVATIONCURVATURE} reads \begin{eqnarray}\label{RBCFIND}
g^{i \overline{j}} g^{p \overline{q}} R_{ \alpha \overline{\beta}\gamma \overline{\delta}}^h f_i^{\alpha} \overline{f_{j}^{\beta}} f_p^{\gamma} \overline{f_{q}^{\delta}} &=& \sum_{\alpha, \gamma} R^h_{\alpha \overline{\alpha} \gamma \overline{\gamma}} \lambda_{\alpha}^2 \lambda_{\gamma}^2.
\end{eqnarray}

This motivated Yang--Zheng \cite{YangZhengRBC} to consider the following (c.f., \cite{LeeStreets}): The \textit{real bisectional curvature} $\text{RBC}_{g}$ of a Hermitian metric $g$ is the function $\text{RBC}_{g} : \mathcal{F}_M \times \mathbb{R}^n \backslash \{ 0 \} \longrightarrow \mathbb{R}$,  $$\text{RBC}_{g}(\vartheta, \lambda) := \frac{1}{| \lambda |^2} \sum_{\alpha, \gamma} R_{\alpha \overline{\alpha} \gamma \overline{\gamma}} \lambda_{\alpha} \lambda_{\gamma}.$$ Here, $\mathcal{F}_M$ denotes the unitary frame bundle, $R_{\alpha \overline{\beta} \gamma \overline{\delta}}$ denote the components of the Chern connection with respect to the frame $\vartheta$, and $\lambda = (\lambda_1, ..., \lambda_n) \in \mathbb{R}^n \backslash \{ 0 \}$.

If the metric is K\"ahler, the real bisectional curvature is comparable to the holomorphic sectional curvature. For a general Hermitian metric, however, the real bisectional curvature strictly dominates the holomorphic sectional curvature.  It is not strong enough, however, to dominate the Ricci curvatures (see \cite{YangZhengRBC}).

It was observed in \cite{BroderSBC} that the real bisectional curvature does not give suitable control in the Aubin--Yau second-order estimate. Indeed, if $f:(M, g) \longrightarrow (N, h)$ is biholomorphic onto its image, taking the trace of \eqref{BochnerMotivation}, with respect to $h$, the corresponding curvature term we want to understand is \begin{eqnarray}\label{AYCONF}
h^{\gamma \overline{\delta}} g^{i \overline{q}} g^{p \overline{j}} R^g_{k \overline{\ell} p \overline{q}} h_{\alpha \overline{\beta}} f_i^{\alpha} \overline{f_j^{\beta}} (f^{-1})_{\gamma}^k \overline{(f^{-1})_{\delta}^{\ell}}.
\end{eqnarray}

Again, choose coordinates at $p$ and $f(p)$ such that $g_{i \overline{j}} =\delta_{ij}$, $h_{\alpha \overline{\beta}} = \delta_{\alpha\beta}$, and $f_i^{\alpha} = \lambda_i \delta_i^{\alpha}$. Then \eqref{AYCONF} reads \begin{eqnarray}\label{AYFIND}
\sum_{i,k} R^g_{i \overline{i} k \overline{k}} \lambda_i^2 \lambda_k^{-2}.
\end{eqnarray}

This is not controlled by the real bisectional curvature (what is the vector here?). 

To understand both \eqref{RBCFIND} and \eqref{AYFIND}, introduce the matrix $\mathcal{R} \in \mathbb{R}^{n \times n}$ with entries $\mathcal{R}_{\alpha \gamma} : = R_{\alpha \overline{\alpha} \gamma \overline{\gamma}}$. If $v = (v_1, ..., v_n) \in \mathbb{R}^n \backslash \{ 0 \}$, then the real bisectional curvature can be written as the Rayleigh quotient \begin{eqnarray*}
\text{RBC}_g(v) & := & \frac{v^t \mathcal{R} v}{v^t v}.
\end{eqnarray*}

From the scale invariance of the Rayleigh quotient, we may equivalently assume that $v \in \mathbb{S}^{n-1} \subset \mathbb{R}^n$. It is then immediate that the maximum of $\text{RBC}_g$ is always achieved, and the maximum occurs precisely when $v$ is the eigenvector of $\mathcal{R}$ corresponding to the largest eigenvalue.

Since the vector $v$ arising in \eqref{RBCFIND} is given by the vector of principal values $$\lambda_1 \geq \lambda_2 \geq \cdots \geq  \lambda_r \geq \lambda_{r+1}  = \cdots = 0,$$ of $\partial f$, where $r = \text{rank}(\partial f)$, it is clear that we need only control the Rayleigh quotient over the cone \begin{eqnarray*}
\Gamma & : = & \{ x \in \mathbb{R}_+^n : x_1 \geq \cdots \geq x_n \geq 0 \}.
\end{eqnarray*}

For the Aubin--Yau inequality, we want to estimate (from below), the quantity \eqref{AYFIND}. Let $\Gamma_+ :=\{(x_1, ..., x_n) \in \mathbb{R}_+^n : x_1 \geq x_2 \geq \cdots \geq x_n > 0 \}$ denote the cone of ordered positive $n$-tuples. For a vector $v \in \Gamma_+$, we denote by $u_v := v_{\circ}^{-1}$ the vector which inverts $v$ with respect to the Hadamard product. That is, if $v = (v_1, ..., v_n) \in \Gamma_+$, then $u_v = (v_1^{-1}, ..., v_n^{-1})$.

Then a bound on \eqref{AYFIND} translates to a bound on the generalized Rayleigh quotient $$\frac{u_v^t \mathcal{R} v}{| u_v |  |v |}.$$

To distinguish this from the real bisectional curvature, we define the following:

\subsection{Definition}\label{SBCDefinition}
Let $(M, g)$ be a Hermitian manifold. Define the \textit{first Schwarz bisectional curvature} \begin{eqnarray*}
\text{SBC}_{g}^{(1)} : \mathcal{F}_M \times \Gamma_+ \to \mathbb{R}, \hspace{0.5cm} \text{SBC}_{g}^{(1)}(\vartheta,v) \ : = \ \frac{u_v^t \mathcal{R} v}{| u_v |  |v |},
\end{eqnarray*}

and the \textit{second Schwarz bisectional curvature} \begin{eqnarray*}
\text{SBC}_{g}^{(2)} : \mathcal{F}_M \times \Gamma \to \mathbb{R}, \hspace*{0.5cm} \text{SBC}_{g}^{(2)}(\vartheta, v) \ : = \ \frac{v^t \mathcal{R} v}{ | v |^2},
\end{eqnarray*}

where $\mathcal{R}$ is the matrix with entries $\mathcal{R}_{\alpha \gamma} : = R_{\alpha \overline{\alpha} \gamma \overline{\gamma}}$ with respect to the frame $\vartheta$.

\hfill

\iffalse
\subsection{Definition}
Let $A \in \mathbb{R}^{n \times n}$ be a symmetric matrix. We say that $A$ is \textit{copositive} (respectively, \textit{conegative}) if $v^t A v \geq 0$ (respectively, $v^t A v \leq 0$) for all $v \in \mathbb{R}_+^n$. We say that $A$ is \textit{strictly copositive} (respectively, \textit{strictly conegative}) if $v^t A v > 0$ (respectively, $v^t A v < 0$ ) for all $v \in \mathbb{R}_+^n$. \footnote{Here, $\mathbb{R}^n_+ : = \{ (x_1, ..., x_n) \in \mathbb{R}^n : x_k > 0, \ \forall k \}$ denotes the positive orthant. Clearly, copositive matrices generalize positive semi-definite matrices. } Let $\Gamma \subset \mathbb{R}_+^n$ be a cone. A symmetric matrix $A \in \mathbb{R}^{n \times n}$ is said to be $\Gamma$--copositive if $v^t A v \geq 0$ for all $v \in \Gamma$.

\subsection{Lemma}
Let $(M, g)$ be a Hermitian manifold. \begin{itemize}
\item[(i)] The first Schwarz bisectional curvature $\text{SBC}_{g}^{(1)}$ is semi-positive (respectively, positive) if, for all frames, the matrix $\mathcal{Q}$ is $\Gamma$--copositive (respectively, strictly $\Gamma$--copositive).
\item[(ii)] The second Schwarz bisectional curvature $\text{SBC}_{g}^{(2)}$ is semi-negative (respectively, negative) if, for all frames, the matrix $\mathcal{Q}$ is $\Gamma$--conegative (respectively, strictly $\Gamma$--conegative).
\end{itemize}

\fi

We now have the tools to give a unified treatment of (and extend) the various forms of the Schwarz lemmas: 

\subsection{Theorem}\label{HermitianCL}
(Hermitian Chern--Lu \cite{BroderSBC,YangZhengRBC}). Let $f : (M^n, g) \longrightarrow (N, h)$ be a holomorphic map of rank $r$ between Hermitian manifolds. Assume $\text{Ric}^{(2)}_{g} \geq - C_1 g - C_2 f^{\ast} h$, and $\text{SBC}^{(2)}_{g} \leq - \kappa$ for some constants such that $\kappa - C_2 > 0$. If $M$ is compact, then \begin{eqnarray*}
| \partial f |^2 & \leq & \frac{C_1  r}{\kappa - C_2}.
\end{eqnarray*}

\subsection{Theorem}
(Hermitian Aubin--Yau \cite{BroderSBC}). Let $f : (M^n, g) \longrightarrow (N, h)$ be a holomorphic map, which is biholomorphic onto its image. Assume $\text{SBC}_{g}^{(1)} \geq - \kappa$ and $\text{Ric}_{h}^{(2)} \leq - C_1 h + C_2 (f^{-1})^{\ast} g$ for constants $\kappa, C_1,C_2$, with $C_1 >0$. Then, if $M$ is compact,  \begin{eqnarray*}
| \partial f |^2 & \leq & \frac{nC_2 + \kappa}{C_1}.
\end{eqnarray*}

\hfill

The Schwarz bisectional curvatures yield an interesting comparison between the Chern--Lu and Aubin--Yau inequalities. Indeed, the Chern--Lu inequality requires an upper bound on $\text{SBC}^{(2)}$, a Rayleigh quotient, which is well-known to give a variational characterization of the eigenvalues. The Aubin--Yau inequality requires a lower bound on $\text{SBC}_{\omega}^{(1)}$,  a generalized Rayleigh quotient, which is known to give a variational characterization of the singular values. Therefore, at least philosophically,  the Chern--Lu inequality is to the Aubin--Yau inequality what the eigenvalue decomposition is to the singular value decomposition.

\subsection{The Chen--Cheng--Lu Schwarz Lemma}
So far, the underlying idea is that a Schwarz lemma is what results from applying a Bochner formula to some elementary symmetric function associated with some holomorphic map. In particular, if one has a Bochner formula, one should get a Schwarz lemma.\footnote{This can be thought of as a `dual' version of Bloch's principle.} The following Schwarz lemma of Chen--Cheng--Lu \cite{ChenChengLu} does not seem to fit into the Bochner paradigm:

\subsection{Theorem}
Let $f : (M ,g) \longrightarrow (N, h)$ be a holomorphic map from a complete K\"ahler manifold $(M,g)$ with $\text{HSC}_{g} \geq \kappa_1$, into a Hermitian manifold $(N, h)$ with $\text{HSC}_{h} \leq \kappa_2 < 0$. If the sectional curvature of $g$ is bounded from below, then \begin{eqnarray*}
| \partial f |^2 & \leq & \frac{\kappa_1}{\kappa_2}.
\end{eqnarray*}

The lower bound on the sectional curvature is, of course, automatic if $M$ is compact, and appears because of the use of a comparison theorem.  Further results in this direction were obtained by Ni \cite{Ni}.  In particular, in \cite[Theorem 1.4]{Ni}, the sectional curvature lower bound is replaced by a lower bound on the bisectional curvature.

Since the Bochner formula will typically yield a Ricci term, the Chen--Cheng--Lu Schwarz lemma is, in a sense, the outlier. If one places holomorphic maps as the center of focus, however, we see that a lower bound on the holomorphic sectional appears in the source of the Aubin--Yau estimate, while an upper bound on the holomorphic sectional curvature appears on the target in the Chern--Lu inequality. This idea leads to the following $8$ real dimensional family of Schwarz lemmas, which includes the Chen--Cheng--Lu Schwarz lemma as a special case:

\subsection{Theorem}
Let $f : (M^n, g) \longrightarrow (N, h)$ be a holomorphic map of rank $r$ between Hermitian manifolds with $\text{SBC}^{(1)}_{g} \geq - \kappa_1$ and $\text{SBC}_{h}^{(2)} \leq - \kappa_2$, for some constants $\kappa_1, \kappa_2 \geq 0$. Assume there is a Hermitian metric $\mu$ on $M$ such that, for constants $C_1, C_2, C_3, C_4 \in \mathbb{R}$, with $C_3>0$, $\kappa_2 - C_2 >0$, we have \begin{eqnarray*}
- C_1 \mu - C_2 f^{\ast} h \ \leq \ \text{Ric}^{(2)}_{\mu} \ \leq \ - C_3 \mu + C_4 g.
\end{eqnarray*}

Then, if $M$ is compact, \begin{eqnarray*}
| \partial f|^2 & \leq & \frac{C_1 n r(\kappa_1 + C_4)}{C_3(\kappa_2 - C_2)}.
\end{eqnarray*}

In particular, the Chen--Cheng--Lu Schwarz lemma fits within the Bochner paradigm.

\subsection*{The Schwarz lemma in other contexts}
Although we have focused exclusively on the role of the Schwarz lemma in complex differential geometry,  the Schwarz lemma appears in many other branches of mathematics: In symplectic geometry,  the Schwarz lemma for pseudo-holomorphic curves (known as \textit{Gromov's Schwarz lemma}) plays an important role (see, e.g., \cite{Gromov, Muller}). The Beardon--Stephenson discrete Schwarz--Pick lemma \cite{BeardonStephenson} for circle packing was used by He--Schramm \cite{HeSchramm} to give a new proof of Thurston's approach to the Riemann mapping theorem by circle packing, previously proven by Rodin--Sullivan \cite{RodinSullivan} (see also \cite{RodinSchwarzLemma}). Let us also mention Tosatti's Schwarz lemma in the almost Hermitian case \cite{TosattiSchwarzLemma} and Dong--Ren--Yu's Schwarz lemmas in the pseudo-Hermitian case \cite{DongRenYu}.  In \cite{Chen}, Chen uses the Schwarz lemma in the Sasakian setting, obtaining an interesting Sasakian extension of the Wu--Yau theorem.

\subsection*{Concluding Remarks}
In the same manner that the Bloch principle predicts the Schwarz lemma from a vanishing or rigidity theorem, we have proposed that a Bochner principle will give rise to a Schwarz lemma. Given the number of possible Bochner formulas that exist, one would suspect an abundance of Schwarz lemmas to appear.

\subsection*{Acknowledgements}
I would like to express my gratitude to my advisors Ben Andrews and Gang Tian for many valuable discussions on the Schwarz lemma and for their constant support and encouragement. In particular,  I'm grateful to Tian for encouraging me to pursue the Schwarz lemma,  striving for a better feeling for these calculations.  I owe a tremendous amount to the innumerable conversations and communications with Yanir Rubinstein and Fangyang Zheng, which have had an important influence on this work.

%\end{multicols}

\hfill

\scshape{Mathematical Sciences Institute, Australian National University, Acton, ACT 2601, Australia}

\scshape{BICMR, Peking University, Beijing, 100871, People's republic of china}

\textit{E-mail address}: \texttt{kyle.broder{@}anu.edu.au}

\end{document}